\documentclass[12pt, leqno]{article}
\usepackage{amsfonts}
\usepackage{amsmath}
\usepackage{amssymb}  
\renewcommand{\baselinestretch}{1}
\begin{document}

\begin{center}{\large\bf BOUNDS FOR ODD $k$-PERFECT NUMBERS}
\vskip1cm
\textsc{Shi-Chao Chen\qquad Hao Luo}
\end{center}

\vskip .5 cm

\parskip = 0.2cm
\topmargin=-0.1cm \oddsidemargin=0cm \evensidemargin=0cm
\renewcommand\baselinestretch{1.0}

\noindent \textsc{Abstract.} Let $k\ge2$ be an integer. A natural
number $n$ is called $k$-perfect if $\sigma(n)=kn.$ For any integer
$r\ge1$ we prove that the number of odd $k$-perfect numbers with at
most $r$ distinct prime factors is bounded by $k4^{r^3}$.\vskip1cm

\begin{center}{ 1. \textsc{Introduction}}\end{center}

Let $\sigma(n)$ be the sum of positive divisors of a natural number
$n$. For a rational number $k>1$, if $\sigma(n)=kn$, then $n$ is
called {\em multiperfect} (or {\em$k$-perfect}). In the special case
when $k=2$, $n$ is called a {\em perfect number}. No odd $k$-perfect
numbers are known for any integer $k\ge2$.

In 1913, Dickson [3] proved that for any natural number $r$, there
are only finitely many odd perfect numbers $n$ with $\omega(n)\le
r$, where $\omega(n)$ is the number of distinct prime factors of the
positive integer $n$.  Pomerance [8] gave an explicit upper bound
$n\le (4r)^{(4r)^{2^{r^2}}}$ in 1977.  Heath-Brown [4] later
improved the bound to $n<4^{4^r}$ and Cook [2] refined this to
$n<{(195)^{\frac{4^r}{7}}}$. In 2003, Nielsen [5] improved the bound
further and proved that for any integer $k\ge2$ if $n$ is an odd
$k$-perfect number with $r$ distinct prime factors then
$$n\le2^{4^r}.\leqno{(1)}$$

Recently, Pollack [7] modified Wirsing's method [9] and bounded the
number of such $n$. He showed that for each positive integer $r$,
the number of odd perfect numbers $n$ with $\omega(n)\le r$ is
bounded by $4^{r^2}$.

In this note we will generalize Pollack's result to odd $k$-perfect
numbers. We have the following

\noindent{\bf Theorem 1.} {\it Let $k\ge2$ be an integer. Then for
any integer $r\ge1$, the number of odd $k$-perfect numbers $n$  with
$\omega(n)\le r$ is bounded by $k4^{r^3}$.} \vskip.5cm

\begin{flushleft}
\medskip
\rule{2.5cm}{0.2mm}
Supported by the Natural Science Foundation of China (Grant
11026080) and the Natural Science Foundation of Education Department
of Henan Province (Grant 2009A110001).
\end{flushleft}

\noindent \begin{center}\textsc{ 2. Proofs}\end{center} \vskip.5cm

If $n_1$ is $k_1$-perfect, $n_2$ is $k_2$-perfect and $(n_1,n_2)=1$,
then $n_1n_2$ is $k_1k_2$-perfect. For this reason, we give the
following

\noindent{\bf Definition 2.} {\it A multiperfect number $n$ is
called primitive if for any $d|n,1<d<n,(d,\frac{n}{d})=1$, we have}
$$d\nmid\sigma(d).$$\vskip.3cm

For example, if $n$ is an odd perfect number, then $n$ is primitive.
The reason is that if there is a divisor $d$ of an odd perfect
number $n$ with $1<d<n, d\mid\sigma(d)$, then
$\frac{\sigma(d)}{d}\ge2.$ Therefore
$$2=\frac{\sigma(n)}{n}
=\sum_{m|n}\frac{m}{n} =\sum_{m|n}\frac{1}{m}
>\sum_{m|d}\frac{1}{m}
=\frac{\sigma(d)}{d}\ge2,$$  which is absurd.

\noindent{\bf Lemma 3.} {\it Let $x\ge1$ and $\alpha>1$ be a
positive rational number. Let $I$ be the number odd primitive
$\alpha$-perfect numbers $n\le x$ with $\omega(n)\le r$. Then
$$I\le1.31\frac{\alpha}{\alpha-1}(\log x)^r.$$If $\alpha$
is an integer, then}
$$I\le0.05(\log x)^r.$$

\noindent{\it Proof.} The proof essentially is in sprite of
Pollack's work [7] and is a modification of Wirsing's method [9].
Let $n\le x$ be an odd primitive $\alpha$-perfect number and
$\omega(n)= s\le r$. Denote by $\nu_p(n)$ the highest power of prime
$p$ dividing $n$. Suppose $p_1$ is the smallest positive prime
factor of $n$ and set $e_1=\nu_{p_1}(n)$. Since $n$ is primitive, we
have
$$\frac{n}{p_1^{e_1}}\nmid\sigma\left(\frac{n}{p_1^{e_1}}\right).$$
It follows that there exist  prime $p_2|\frac{n}{p_1^{e_1}}$ such
that
$$\nu_{p_2}\left(\frac{n}{p_1^{e_1}}\right)
>\nu_{p_2}\left(\sigma\left(\frac{n}{p_1^{e_1}}\right)\right).$$
We suppose that $p_2$ is the smallest such prime and write
$e_2=\nu_{p_2}(n)$. Replacing $\frac{n}{p_1^{e_1}}$ by
$\frac{n}{p_1^{e_1}p_2^{e_2}}$ and iterating the argument above, we
can determine $p_3$. Write $e_3=\nu_{p_3}(n)$. Continuing this
procedure, we can obtain primes $p_i$ and exponents
$e_i=\nu_{p_i}(n), i=4\cdots,s$. Hence $n$ has standard
factorization
$$n=p_1^{e_1}p_2^{e_2}\cdots p_s^{e_s}.$$
We need count the number of possibilities for such primes $p_i$ and
exponents $e_j$. The algorithm shows that $p_2$ is determined by
$p_1$ and $e_1$, $p_3$ is determined by $p_1,p_2$ and $e_1,e_2$,
$p_i$ is determined by $p_1,\cdots, p_{i-1}$ and
$e_1,\cdots,e_{i-1}$. Therefore it is sufficient to count the number
of possibilities of $p_1$ and $e_1,e_2,\cdots,e_s$. We have
\begin{align*}\alpha&=\frac{\sigma(n)}{n}
=\prod_{p^e||n}\frac{p^{e+1}-1}{p^e(p-1)}
=\prod_{p^e||n}\frac{p^{e+1}-1}{p^e(p-1)}
=\prod_{p^e||n}\left(1+\frac{p^{e}-1}{p^e(p-1)}\right)\\
&<\prod_{p|n}\left(1+\frac{1}{(p-1)}\right)
=\prod_{p|n}\left(1-\frac{1}{p}\right)^{-1}
\le\left(1-\sum_{p|n}\frac{1}{p}\right)^{-1}
<\left(1-\frac{s}{p_1}\right)^{-1}.\end{align*} Thus
$$p_1\le\frac{\alpha s}{\alpha-1}.$$
Since $p_i^{e_i}||n, n\le x$, we have
$$e_i\le\frac{\log x}{\log p_i}.$$
The number of possibilities for the sequence $p_1,e_1,\cdots,e_s$ is
bounded by
$$\pi\left(\frac{\alpha s}{\alpha-1}\right)\prod_{i=1}^s\frac{\log x}{\log
p_i},$$where $\pi(\frac{\alpha s}{\alpha-1})$ is the number of
primes not exceeding $\frac{\alpha s}{\alpha-1}$. Since $1\le
s=\omega(n)\le r$, it follows that
\begin{align*}I&\le r\cdot\pi\left(\frac{\alpha
s}{\alpha-1}\right)\prod_{i=1}^s\frac{\log x}{\log p_i}\\
&\le\frac{\frac{\alpha }{\alpha-1}r^2}{2\log p_1\log p_2\cdots\log
p_r}(\log x)^r\\
&\le\frac{\frac{\alpha }{\alpha-1}r^2}{2\log q_1\log q_2\cdots\log
q_r}(\log x)^r,\tag{2}\end{align*} where $q_i$ is the $i$-th odd
prime, $q_1=3, q_2=5,\cdots.$ For convenience, we denote by
$$f(r):=\frac{r^2}{2\log q_1\log
q_2\cdots\log q_r}.$$ By simple calculation, $f(r)$ is a decreasing
function of $r$ for $r\ge3$. The maximal value of $f(r)$ is
$$ f(3)=\frac{9}{2\log3\log5\log7}<1.31. \leqno{(3)}$$
If $\alpha=2$, then Nielsen [6] showed that $\omega(n)\ge9$ for any
odd perfect $n$. If $\alpha\ge3$ is an integer and $n$ is an odd
$\alpha$-perfect number, then Cohen and Hagis [1] proved that
$\omega(n)\ge11$. It follows that for any integer $\alpha\ge2$,
$$ f(r)\le f(9)=\frac{9}{2\log 3\log 5\cdots\log 29}< 0.022. \leqno{(4)}$$
Lemma 3 follows from (2), (3) and (4).\quad$\Box$

\noindent{\bf Lemma 4.} {\it Let $x\ge1, r\ge1$ and integer $k\ge2$.
The number of odd $k$-perfect $n\le x$ with $\omega(n)\le r$ is
bounded by $k(\log x)^{\frac{r^2+8r}{9}}$.}

\noindent{\it Proof.} Suppose that $\sigma(n)=kn$. Let $d_1$ be the
smallest positive divisors of $n$ with
$1<d_1<n,(d_1,\frac{n}{d_1})=1$ and $d_1|\sigma(d_1)$. Then it is
easy to see that $d_1$ is an odd primitive multiperfect number. We
write $\sigma(d_1)=k_1d_1$ for some integer $k_1$. Similarly, let
$d_2$ be the smallest positive divisors of $\frac{n}{d_1}$ with
$1<d_2<\frac{n}{d_1},(d_2,\frac{n}{d_1d_2})=1$ and
$d_2|\sigma(d_2)$. Then  $d_2$ is an odd primitive multiperfect
number. We write $\sigma(d_2)=k_2d_2$ for some integer $k_2$.
Continuing this algorithm, we can find divisors $d_i$ of $n,
i=1,\cdots, j,$ and integers $k_i,i=1,2,\cdots j$ such that
$$d_i|\frac{n}{d_1\cdots d_{i-1}},\quad \left(d_i,\frac{n}{d_1\cdots
d_{i-1}d_i}\right)\ and\ \sigma(d_i)=k_id_i\ for\ some\ integer\
k_i\ge2.$$ We assume that the procedure stops at the $j+1$-th step
when $\frac{n}{d_1d_2\cdots d_j}=1$ or $\frac{n}{d_1d_2\cdots d_j}$
is primitive and $\frac{n}{d_1d_2\cdots
d_j}\nmid\sigma(\frac{n}{d_1d_2\cdots d_j})$. Definite $d_{j+1}$ by
$$d_{j+1}=\frac{n}{d_1d_2\cdots d_j}.$$ Then we have
$$n=d_1d_2\cdots d_jd_{j+1}.\leqno{(5)}$$
If $d_{j+1}\neq1$, then
\begin{align*}kn&=\sigma(n)\\
&=\sigma(d_1d_2\cdots d_{j+1})\\
&=\sigma(d_1)\sigma(d_2)\cdots\sigma(d_{j+1})\\
&=k_1d_1k_2d_2\cdots k_jd_j\sigma(d_{j+1}).\end{align*} Therefore
$$\sigma(d_{j+1})=\frac{k}{k_1k_2\cdots k_j} d_{j+1}.$$
It follows that $d_{j+1}$ is  $\frac{k}{k_1k_2\cdots k_j}$-perfect
and ${k_1k_2\cdots k_j}\nmid k$. Since $k_1,\cdots,k_s$ are
integers, we have  $$k_1k_2\cdots k_j\le k-1.$$ In view of Lemma 3,
the number of such $d_{j+1}$ not exceeding $x$ is bounded by
$$1.31\frac{\frac{k}{k_1k_2\cdots k_j}}{\frac{k}{k_1k_2\cdots k_j}-1}(\log x)^r
=1.31\frac{k}{k-k_1k_2\cdots k_j}(\log x)^r\le1.31k(\log x)^r.$$ By
the minimality of $d_1,\cdots,d_j$, one can see that all $d_1,
\cdots, d_j$ are primitive. Nielsen [6], Cohen and Hagis's results
[1] imply that $\omega(d_i)\ge9, i=1,\cdots,j$. Therefore
$$r\ge\omega(n)=\omega(d_{j+1})+\sum_{i=1}^j\omega(d_i)\ge 1+9j.$$ It
follows that $$j\le\frac{r-1}{9}.$$ By (5) and Lemma 3, the number
of $k$-perfect numbers $n\le x$ with $\omega(n)\le r$ is at most
\begin{align*}(0.05(\log x)^r)^{j}(1.31k(\log x)^r)&\le(0.05(\log
x)^r)^{\frac{r-1}{9}}(1.31k(\log x)^r)\\
&\le \frac{1}{2}k(\log x)^{\frac{r^2+8r}{9}}.\end{align*} If
$d_{j+1}=1$, then $j\le\frac{r}{9}$ and the bound is
$$(0.05(\log x)^r)^{j}\le(\log x)^{\frac{r^2}{9}}\le \frac{1}{2}k(\log x)^{\frac{r^2}{9}}.$$
This complete the proof of Lemma 4.\quad$\Box$

\noindent{\it Proof of Theorem 1.} Let $x=2^{4^r}$. Applying Lemma 4
and Nielson's bound (1), we deduce that the number of odd
$k$-perfect numbers $n$ with $\omega(n)\le r$ is at most
$$k(\log x)^{\frac{r^2+8r}{9}}<k(4^r)^{\frac{r^2+8r}{9}}=k4^{\frac{r^3+8r^2}{9}}\le k4^{r^3}.\quad\Box$$
\vskip.5cm

\begin{center}\textsc{ References}\end{center}

\begin{itemize}

\item[{[1]}] G. L. Cohen and P. Hagis, Jr. {\it Results concerning odd multiperfect
numbers}, Bull. Malaysian Math. Soc. 8(1985), 23-26.

\item[{[2]}] R. J. Cook, {\it Bounds for odd perfect numbers}, In {\it number theory (Ottawa, ON, 1996)},
CRM Proc. Lecture Notes, Vol 19, Amer. Math. Soc., Providence, RI,
1999, 67-71.

\item[{[3]}] L. E. Dickson, {\it Finiteness of the odd perfect and primitive abundant
numbers with $n$ distinct prime factors}, Amer. J. Math. 35(1913),
413-422.

\item[{[4]}] D. R. Heath-Brown, {\it Odd perfect numbers}, Math. Proc. Cambridge Phil.
Soc. 115(1994), 191-196.

\item[{[5]}] P. Nielsen, {\it An upper bound for odd perfect numbers}, Integers:
Electronic J. Comb. Number Theory, 3(2003), A14, 9pp. (electronic).

\item[{[6]}] P. Nielsen, {\it Odd perfect numbers have at least nine distinct prime
factors}, Math. Comp. 76(2007), 2109-2126.

\item[{[7]}] P. Pollack. {\it On Dickson's theorem concerning odd perfect numbers},
Amer. Math. Monthly. 118(2011), 161-164.

\item[{[8]}] C. Pomerance, {\it Multiply perfect numbers, Mersenne primes and
effective computability}, Math. Ann. 226(1977), 195-206.

\item[{[9]}] E. Wirsing, {\it Bemerkung zu der arbeit $\ddot{u}$ber vollkommene Zahlen},
Math. Ann. 137(1959), 316-318.

\end{itemize}
\vskip.5cm

\textsc{Institute of Applied Mathematics, School of Mathematics and
Information Science, Henan University, Kaifeng, 475004, P.R.China}

{\it E-mail address}: Shi-Chao Chen: \texttt{schen@henu.edu.cn}

{\it E-mail address}: Hao Luo:  \texttt{luohao200681@126.com}

\end{document}